\RequirePackage{ifpdf}
\ifpdf 
\documentclass[pdftex]{sigma}
\else
\documentclass{sigma}
\fi

\begin{document}

\allowdisplaybreaks

\renewcommand{\thefootnote}{$\star$}

\renewcommand{\PaperNumber}{068}

\FirstPageHeading

\ShortArticleName{Recurrence Coef\/f\/icients of a New Generalization of the Meixner Polynomials}

\ArticleName{Recurrence Coef\/f\/icients of a New Generalization\\
 of the Meixner Polynomials\footnote{This paper is a
contribution to the Special Issue ``Relationship of Orthogonal Polynomials and Special Functions with Quantum Groups and Integrable Systems''. The
full collection is available at
\href{http://www.emis.de/journals/SIGMA/OPSF.html}{http://www.emis.de/journals/SIGMA/OPSF.html}}}

\Author{Galina FILIPUK~$^\dag$ and Walter VAN ASSCHE~$^\ddag$}

\AuthorNameForHeading{G.~Filipuk and W.~Van Assche}

\Address{$^\dag$~Faculty of Mathematics, Informatics and Mechanics,
University of Warsaw,\\
\hphantom{$^\dag$}~Banacha 2, Warsaw, 02-097, Poland}
\EmailD{\href{mailto:filipuk@mimuw.edu.pl}{filipuk@mimuw.edu.pl}}
\URLaddressD{\url{http://www.mimuw.edu.pl/~filipuk/}}

\Address{$^\ddag$~Department of Mathematics, Katholieke Universiteit Leuven,\\
\hphantom{$^\ddag$}~Celestijnenlaan 200B box 2400, BE-3001 Leuven, Belgium}
\EmailD{\href{mailto:walter@wis.kuleuven.be}{walter@wis.kuleuven.be}}
\URLaddressD{\url{http://wis.kuleuven.be/analyse/walter/}}

\ArticleDates{Received April 18, 2011, in f\/inal form July 07, 2011;  Published online July 13, 2011}

\Abstract{We investigate new generalizations of the Meixner polynomials on
the lattice $\mathbb{N}$, on the shifted lattice
$\mathbb{N}+1-\beta$ and on the bi-lattice $\mathbb{N}\cup
(\mathbb{N}+1-\beta)$. We show that the coef\/f\/icients of the
three-term recurrence relation for the orthogonal polynomials are
related to  the solutions of the  f\/ifth Painlev\'e equation {\rm P}$_{\textup V}$.
 Initial conditions for dif\/ferent lattices can be transformed to the classical solutions
of {\rm P}$_{\textup V}$ with  special values of the parameters. We also study
one property of the B\"acklund transformation of {\rm P}$_{\textup V}$.}


\Classification{34M55; 33E17; 33C47; 42C05; 64Q30}

\renewcommand{\thefootnote}{\arabic{footnote}}
\setcounter{footnote}{0}

\section{Introduction}

The recurrence coef\/f\/icients of orthogonal polynomials for
semi-classical weights are  often related to the  Painlev\'e type
equations (e.g., \cite{Magnus, GFWalterLaguerre, GFWalterCharlier} and see
also the overview in \cite{LiesGFWalter}). In this paper we study a
new generalization of the Meixner weight. The recurrence
coef\/f\/icients
  of the corresponding
orthogonal polynomials can be viewed as functions of one of the
parameters. We show that the recurrence coef\/f\/icients are related
to solutions of the  f\/ifth Painlev\'e equation. Another
generalization of the Meixner weight is presented in
\cite{LiesGFWalter}.

The paper is organized as follows. In the introduction we shall
f\/irst review orthogonal polynomials  for the generalized Meixner
weight on dif\/ferent lattices and their main properties fol\-lowing~\cite{bilattice}. Next we shall brief\/ly recall the f\/ifth
Painlev\'e equation and its B\"acklund transformation. Further, by
using the Toda system, we show that the recurrence coef\/f\/icients
can be expressed in terms of solutions of the f\/ifth Painlev\'e
equation. Finally  we study initial conditions of  the recurrence
coef\/f\/icients for dif\/ferent lattices and describe one property of
the B\"acklund transformation of {\rm P}$_{\textup V}$.

\subsection{Orthogonal polynomials for the generalized Meixner weight}

One of the most important properties of orthogonal polynomials is
the three-term recurrence relation. Let us consider a sequence
$(p_n)_{n \in \mathbb{N}}$ of orthonormal polynomials for the weight
$w$ on the lattice $\mathbb{N}=\{0,1,2,3,\ldots\}$
\begin{gather*}\label{orthonormality}
\sum_{k=0}^{\infty} p_n(k)p_k(k)w(k)=\delta_{n,k},
\end{gather*}
where $\delta_{n,k}$ is the Kronecker delta. This relation takes
the following form:
\begin{gather}\label{3 term orthonormal}
 x p_n(x)=a_{n+1}p_{n+1}(x)+b_np_n(x)+a_np_{n-1}(x).
\end{gather}
 For the   monic
polynomials $P_n$ related to orthonormal polynomials
$p_n(x)=\gamma_n x^n+\cdots $ with
\[\frac{1}{\gamma_n^2}=\sum_{k=0}^{\infty}P_n^2(k)w(k), \] the
recurrence relation is given by \[x
P_n(x)=P_{n+1}(x)+b_nP_n(x)+a_n^2 P_{n-1}(x).\]

The classical Meixner polynomials (\cite[Chapter VI]{Chihara},
\cite{bilattice}) are orthogonal on the lattice $\mathbb{N}$ with
respect to the
 negative binomial (or Pascal) distribution: %
\[
 \sum_{k=0}^{\infty}M_n(k;\beta,c)M_m(k;\beta,c)\frac{(\beta)_k c^k}{k!}=\frac{c^{-n}n!}{(\beta)_n(1-c)^{\beta}}\,\delta_{n,m},\qquad \beta>0,\quad 0<c<1.
\]
Here the Pochhammer symbol is def\/ined by \[(x)_n=\frac{\Gamma(x+n)}{\Gamma(x)}=x(x+1)\cdots(x+n-1).\]
The weight $w_k=w(k)=(\beta)_k c^k/k!$ satisf\/ies the Pearson
equation \[\nabla [(\beta+x)
w(x)]=\left(\beta+x-\frac{x}{c}\right)w(x) ,\] where $\nabla$
is the backward dif\/ference operator \[\nabla f(x)=f(x)-f(x-1).\]
Here the  function
$w(x)=\Gamma(\beta+x)c^x/(\Gamma(\beta)\Gamma(x+1))$ gives the
weights $w_k=w(k).$ The
 Pearson equation for the Meixner polynomials is, hence, of the form
\begin{gather}\label{Pearson}\nabla[\sigma(x)w(x)]=\tau(x)w(x)
\end{gather}
with $\sigma(x)=\beta+x$ and $\tau$
 is a polynomial of degree 1.

  Recall that the classical
 orthogonal polynomials are characterized by the Pearson equation~ (\ref{Pearson})
 with $\sigma$ a polynomial of degree at most~2 and~$\tau$
 a polynomial of degree~1. Note that in (\ref{Pearson}) the
 operator
 $\nabla$ is used for orthogonal polynomials on the lattice and
   it is replaced by dif\/ferentiation in case of orthogonal
  polynomials on an interval of the real line. The Pearson equation plays an
 important role for classical orthogonal polynomials since it
 allows to f\/ind many useful properties of these polynomials.
 It is known  that the recurrence coef\/f\/icients of the Meixner polynomials are given explicitly by
\[a_n^2=\frac{n(n+\beta-1)c}{(1-c)^2},\qquad b_n=\frac{n+(n+\beta)c}{1-c},\qquad n \in\mathbb{N}.\]

 The Meixner  weight can be generalized \cite{bilattice}. One can
 use the weight function \[w(x)=\frac{\Gamma(\beta)\Gamma(\gamma+x)c^x}{\Gamma(\gamma)\Gamma(\beta+x)\Gamma(x+1)}
\] which gives the weight \begin{gather}\label{weight}
 w_k=w(k)=\frac{(\gamma)_k c^k}{(\beta)_k
 k!},\qquad c,\beta,\gamma>0,\end{gather} on the lattice $\mathbb{N}.$
  The orthonormal   polynomials  $(p_n)_{n\in\mathbb{N}}$ for weight (\ref{weight})
 satisfy
\begin{gather}\label{orthogonality1}\sum_{k=0}^{\infty}p_n(k)p_m(k)w_k=0,\qquad n\neq
m.\end{gather}

The special case $\beta=1$ was studied in \cite{LiesGFWalter}. The
case $\beta=\gamma$ gives the well-known Charlier weight.   The
case $\gamma=1$ corresponds to the  classical Charlier weight on
the shifted lattice $\mathbb{N}+1-\beta$.

\begin{theorem}[\protect{\cite[Theorem 3.1]{bilattice}}]\label{thm:
discrete} The recurrence coefficients in the three-term recurrence
relation \eqref{3 term orthonormal} for the orthonormal
polynomials defined by \eqref{orthogonality1}, with respect to
weight \eqref{weight} on the lattice~$\mathbb{N}$, satisfy
\[a_n^2=nc-(\gamma-1)u_n,\qquad b_n=n+\gamma-\beta+c-(\gamma-1)v_n/c,\]
where
\begin{gather}\label{discrete1}
(u_n+v_n)(u_{n+1}+v_n)=\frac{\gamma-1}{c^2}v_n(v_n-c)\left(v_n-c\frac{\gamma-\beta}{\gamma-1}\right),
\\ \label{discrete2}
  (u_n+v_n)(u_n+v_{n-1})=\frac{u_n}{u_n-c
n/(\gamma-1)}(u_n+c)\left(u_n+c\frac{\gamma-\beta}{\gamma-1}\right),
\end{gather}
with initial conditions
\begin{gather}\label{initial}
a_0^2=0,\qquad b_0=\frac{\gamma c}{\beta}
\frac{M(\gamma+1,\beta+1,c)}{M(\gamma,\beta,c)},
\end{gather}
where $M(a,b,z)$ is the confluent hypergeometric function
$_1F_1(a;b;z)$.
\end{theorem}

The system (\ref{discrete1}), (\ref{discrete2})   can be
identif\/ied as a limiting case of an asymmetric discrete Painlev\'e
equation \cite{bilattice}. In this paper we show that it can be obtained from the
B\"acklund transformation of the f\/ifth Painlev\'e equation.
 Furthermore, one can use the
 weight~(\ref{weight}) on the shifted lattice $\mathbb{N}+1-\beta$ and one can
 also combine both lattices to obtain the bi-lattice
 $\mathbb{N}\cup(\mathbb{N}+1-\beta)$. The orthogonality
 measure for the bi-lattice is a linear combination of the
 measures on $\mathbb{N}$ and $\mathbb{N}+1-\beta$.

The weight $w$  in (\ref{weight}) on the  shifted lattice
$\mathbb{N}+1-\beta=\{1-\beta,2-\beta,3-\beta,\ldots\}$ is, up to
a constant factor, equal to the weight on the original lattice
$\mathbb{N}$, with dif\/ferent parameters \cite{bilattice}. Denoting
\[w_{\gamma,\beta,c}(x)=\frac{\Gamma(\beta)}{\Gamma(\gamma)}\frac{\Gamma(\gamma+x)c^x}{\Gamma(x+1)\Gamma(\beta+x)}\]
one has
\begin{gather}\label{weight sh}
w_{\gamma,\beta,c}(k+1-\beta)=c^{1-\beta}\frac{\Gamma(\beta)\Gamma(\gamma+1-\beta)}
{\Gamma(2-\beta)\Gamma(\gamma)}w_{\gamma+1-\beta,2-\beta,c}(k).
\end{gather}
 The corresponding
orthonormal  polynomials $(q_n)_{n\in \mathbb{N}}$  satisfy
\[\sum_{k=0}^{\infty}q_n(k+1-\beta)q_{m}(k+1-\beta)w(k+1-\beta)=0,\qquad n\neq
m. \]  Moreover, these polynomials are equal to the polynomials
$p_n$ shifted in both the variable $x$ and the parameters $\beta$
and  $\gamma$. For the positivity of the weights
$(w(k+1-\beta))_{k \in\mathbb{N}}$ it is necessary to have
$c>0$, $\beta<2$, $\gamma>\beta-1$.  In \cite[Theorem~3.2]{bilattice} it
is shown that the recurrence coef\/f\/icients in the three-term
recurrence relation
\[x q_n(x)=\hat{a}_{n+1}q_{n+1}(x)+ {\hat{b}}_nq_n(x)+\hat{a}_nq_{n-1}(x) \]
satisfy the same system (\ref{discrete1}), (\ref{discrete2}) (with
hats) but with initial conditions
\begin{gather}\label{initial sh}\hat{a}_0^2=0,\qquad {\hat{b}}_0= (1-\beta)\frac{M(\gamma-\beta+1,1-\beta,c)}{M(\gamma-\beta+1,2-\beta,c)}.\end{gather}

Using the orthogonality measure $\mu=\mu_1+\tau \mu_2,$ where
$\tau>0,$ $\mu_1$ is the discrete measure on $\mathbb{N}$ with
weights $w_k=w(k)$ and $\mu_2$ is the discrete measure on
$\mathbb{N}+1-\beta$ with weights $v_k=w(k+1-\beta)$, one can
study orthonormal polynomials $(r_n)_{n\in\mathbb{N}},$
satisfying the three-term recurrence relation
\[x r_n(x)=\tilde{a}_{n+1}r_{n+1}(x)+\tilde{b}_n
r_n(x)+\tilde{a}_nr_{n-1}(x).\] One needs to impose the conditions
$c>0$, $0<\beta<2$, $\gamma>\max(0,\beta-1)$ for the positivity
of the measures.
  The orthogonality relation is
given by
\[\sum_{k=0}^{\infty}r_n(k)r_m(k)w(k)+\tau\sum_{k=0}^{\infty}r_n(k+1-\beta)r_m(k+1-\beta)w(k+1-\beta)=0,\qquad m\neq
n. \] According to \cite[Theorem~3.3]{bilattice} the recurrence
coef\/f\/icients $\tilde{a}_n^2$ and $\tilde{b}_n$ satisfy  system
(\ref{discrete1}), (\ref{discrete2}) (with tilde) but with initial
conditions
\begin{gather}\label{initial
bi}\tilde{a}_0^2=0,\qquad {\tilde{b}}_0= \frac{m_1+\tau
\hat{m}_1}{m_0+\tau \hat{m}_0},
\end{gather}
where
\begin{gather*}
m_0=M(\gamma,\beta,c),\qquad m_1=\frac{\gamma c}{\beta}M(\gamma+1,\beta+1,c),\\
\hat{m}_0=\frac{\Gamma(\beta)\Gamma(\gamma-\beta+1)}{\Gamma(\gamma)\Gamma(2-\beta)}c^{1-\beta}M(\gamma-\beta+1,2-\beta,c),\\
\hat{m}_1=\frac{\Gamma(\beta)\Gamma(\gamma-\beta+1)}{\Gamma(\gamma)\Gamma(1-\beta)}c^{1-\beta } M(\gamma-\beta+1,1-\beta,c).
\end{gather*}

Thus, it is shown in~\cite{bilattice} that the orthogonal
polynomials for the generalized Meixner  weight
 on the
lattice $\mathbb{N}$, on the shifted lattice $\mathbb{N}+1-\beta$
and on the bi-lattice $\mathbb{N}\cup (\mathbb{N}+1-\beta)$ have
recurrence coef\/f\/icients $a_n^2$ and $b_n$ which satisfy the same
nonlinear system of discrete (recurrence) equations  but the
initial conditions are dif\/ferent in each case.

\subsection{The f\/ifth Painlev\'e equation and its B\"acklund transformation}

The Painlev\'e equations possess the so-called Painlev\'e
property: the only movable singularities of the solutions are
poles~\cite{GrLSh}. They are often referred to as nonlinear
special functions and have numerous applications in mathematics
and mathematical physics.

The f\/ifth Painlev\'e equation {\rm P}$_{\textup V}$  is given by
\begin{gather}\label{P5}
y''  =
\left( \frac{1}{2y}+ \frac{1}{y-1}\right)(y')^2
-
 \frac{y'}{t}+ \frac{(y-1)^2}{t^2}\left(A
y +  \frac{B}{y}\right)+ \frac{C
y}{t}+ \frac{D y (y+1)}{y-1},
\end{gather}
where $y=y(t)$ and $A$, $B$, $C$, $D$ are arbitrary complex
parameters.  By using a transformation $y(t)\rightarrow y(k_1 t) $
we can take the value of the parameter $D$ equal  to any non-zero
number. There exists  a B\"acklund transformation between solutions
of the f\/ifth Painlev\'e  equation with $D\neq 0$.

\begin{theorem}[\protect{\cite[Theorem~39.1]{GrLSh}}] \label{thm: BTr2}
If $y=y(t)$ is the solution the fifth Painlev\'e equation
\eqref{P5} with parameters $A$, $B$, $C$, $D$, then the
transformation
\[T_{\varepsilon_1,\varepsilon_2,\varepsilon_3}: \ y\rightarrow
y_1\]  gives another solution $y_1=y_1(t)$ with new values of the
parameters $A_1$, $B_1$, $C_1$, $D_1$, where
\begin{gather*}
y_1=1-\frac{2 d t y}{t  y'-a y^2+(a-b+d t)y+b},\\
 A_1=-\frac{1}{16D}(C+d(1-a-b))^2,\qquad
B_1=\frac{1}{16D}(C-d(1-a-b))^2,\\   C_1=d(b-a),\qquad D_1=D,
\\
 a=\varepsilon_1\sqrt{2A},\qquad b=\varepsilon_2\sqrt{-2B},\qquad d=\varepsilon_3\sqrt{-2D},\qquad
\varepsilon_j^2=1,\qquad j\in\{1,2,3\}.
\end{gather*}
\end{theorem}

 See also \cite{Noumi} for a
further description of the B{\"a}cklund transformations and the
isomorphism of the group of B{\"a}cklund transformations and the
af\/f\/ine Weyl group of $A_3^{(1)}$ type.

\section{Main results}

In this paper we show how to obtain a relation between the recurrence
coef\/f\/icients and the (classical) solutions of the f\/ifth Painlev\'e equation.
The calculations are similar to calculations in~\cite{LiesGFWalter} but are more involved.
 The study of initial conditions of the recurrence coef\/f\/icients for
dif\/ferent lattices is also presented. We can summarize the known
results and our recent f\/indings concerning the (generalized)
Meixner weights as follows.

 The weight $(\beta)_k c^k/(k!)$, $\beta>0$, $0<c<1$, is the classical Meixner
weight and the recurrence coef\/f\/icients are known explicitly. The
weight $(\beta)_k c^k/(k!)^2,$ $\beta>0$, $c>0$, is studied in~\cite{LiesGFWalter} and the recurrence    coef\/f\/icients are related
to classical solutions of  {\rm P}$_{\textup V}$  with parameters
$((\beta-1)^2/2,-(\beta+n)^2/2,2n,-2)$. It is shown in this paper
that the recurrence coef\/f\/icients for the  weight $(\gamma)_k
c^k/(k!(\beta)_k)$,  $c,\beta,\gamma>0,$ are related to the classical
solutions of  {\rm P}$_{\textup V}$  with parameters
$((\gamma-1)^2/2,-(\gamma-\beta+n)^2/2,k_1(\beta+n),-k_1^2/2))$,
$k_1\neq 0$.

\subsection{Relation to the f\/ifth Painlev\'e equation and its B\"acklund transformation}

First we obtain a nonlinear discrete equation for $v_n(c)$. From
equation~(\ref{discrete1}) with $n$ and equation~(\ref{discrete2})
with $n+1$ we eliminate $u_{n+1}$ by computing the resultant.
Next, from the obtained equation and~(\ref{discrete2}) with $n$ we
eliminate~$u_n$. As a result, we obtain a nonlinear discrete
equation for $v_n(c)$ which we denote by
\begin{gather}\label{nonlinear discrete vn}
F(v_{n-1},v_{n},v_{n+1},c)=0.
\end{gather}
The equation was obtained by using
\texttt{Mathematica}\footnote{\url{http://www.wolfram.com}} but it is too long
and too complicated too include here explicitly (all enquiries concerning computations can be sent to the f\/irst author). We shall show
later on that equation (\ref{nonlinear discrete vn}) can in fact
be obtained from the  B\"acklund transformation of the f\/ifth
Painlev\'e equation.

Next we derive the dif\/ferential equation for $v_n$. In
\cite{LiesGFWalter} we have used the Toda system. Since the weight
$w$  in (\ref{weight}) on the  shifted lattice
$\mathbb{N}+1-\beta=\{1-\beta,2-\beta,3-\beta,\ldots\}$ is, up to
a constant factor, equal to the weight on the original lattice
$\mathbb{N}$ with dif\/ferent parameters \cite{bilattice},   it can be shown
\cite{LiesGFWalter, GFWalterCharlier} that the recurrence
coef\/f\/icients $a_n$ and $b_n$ as functions of the parameter $c$
satisfy the Toda system
\begin{gather} \label{Toda}
\begin{split}
& \big(a_n^2\big)':=\frac{d}{dc}\left(a_n^2\right)=\frac{a_n^2}{c}(b_{n}-b_{n-1}), \\
& b_n':=\frac{d}{dc}b_n=\frac{1}{c}\big(a_{n+1}^2-a_n^2\big).
\end{split}
\end{gather}
The same system (\ref{Toda}) holds for the initial lattice  $\mathbb{N}$~\cite{LiesGFWalter}.

Solving  (\ref{discrete1}) for $u_{n+1}$ and  (\ref{discrete2})
for $v_{n-1}$ and substituting into the Toda system (\ref{Toda})
(where we have replaced $a_n^2$ and $b_n$ by their expressions in
terms of $u_n$ and $v_n$ from Theorem  \ref{thm: discrete}), we
get two equations \[
u_n'=R_1(u_n,v_n,c)
\] and
\begin{gather}\label{*}
v_n'=R_2(u_n,v_n,c),
\end{gather}
where the dif\/ferentiation is with respect to $c$. The explicit
expressions for $R_1$ and $R_2$ are again quite complicated but
readily computed in \texttt{Mathematica}. By dif\/ferentiating
equation~(\ref{*}) and substituting the expression for $u_n'$ we
obtain an equation for $v_n''$, $v_n'$, $v_n$, $u_n$. Finally, elimina\-ting~$u_n$ between this equation and~(\ref{*}) gives a nonlinear second
order second degree equation for~$v_n$:
\begin{gather}\label{nonlinear diff vn}
G(v_n'',v_n',v_n,c)=0.
\end{gather}
We have again used \texttt{Mathematica} to compute this long
expression. Now the main dif\/f\/iculty is in identifying this
equation. Since (\ref{discrete1}), (\ref{discrete2}) is similar to
the discrete system in~\cite{LiesGFWalter}, we can try to reduce
equation (\ref{nonlinear diff vn}) to the f\/ifth Painlev\'e
equation. First, we scale the independent variable $c\rightarrow
c/k_1$ and denote $c=k_1 t$, $v(c)=V(t)$. The equation~(\ref{nonlinear diff vn}) becomes $G_1(V_n'',V_n',V_n,t)=0,$ where
the dif\/ferentiation is with respect to $t$. By considering the
Ansatz
\[
V_n(t)=\frac{p_1(t)y'+p_2(t)y^2+p_3(t)y+p_4(t)}{y(y-1)},
\] where
$y=y(t)$ is the solution of the f\/ifth Painlev\'e equation and
$p_j(t)$ are unknown functions, we f\/inally get the following
theorem.

\begin{theorem}\label{thm: P5}
The equation $G_1(V_n'',V_n',V_n,t)=0$ is reduced to the fifth
Painlev\'e equation {\rm P}$_{\textup V}$  by the following transformations:
\begin{enumerate}\itemsep=0pt
\item[$1.$] \begin{gather}\label{case1}
V(t)=\frac{k_1 t(ty'-(1+\beta-2\gamma)y^2+(1+n-k_1 t+\beta-2\gamma
)y-n)}{2(\gamma-1)(y-1)y},
\end{gather}
where $y=y(t)$ satisfies {\rm P}$_{\textup V}$  with
\begin{gather}\label{case1 par}
A=\frac{(\beta-1)^2}{2},\qquad B=-\frac{n^2}{2},\qquad C=k_1(n-\beta+2\gamma),\qquad D=-\frac{k_1^2}{2};\end{gather}
\item[$2.$] \begin{gather}\label{case2}
V(t) = \frac{k_1t(ty'-(\beta-\gamma)y^2+(n-1-k_1 t+\beta
)y+1-n-\gamma)}{2(\gamma-1)(y-1)y},
\end{gather}
where $y=y(t)$ satisfies {\rm P}$_{\textup V}$  with
\begin{gather}\label{case2 par}
A=\frac{(\beta-\gamma)^2}{2},\qquad B=-\frac{(\gamma+n-1)^2}{2},\qquad C=k_1(2+n-\beta),\qquad D=-\frac{k_1^2}{2};\!\!\!\end{gather}
\item[$3.$] \begin{gather*}\label{case3}
V(t)=\frac{k_1t(t
y'+(\gamma-1)y^2+(1+n-k_1t-\beta)y-n+\beta-\gamma)}{2(\gamma-1)(y-1)y},
\end{gather*}
where $y=y(t)$ satisfies {\rm P}$_{\textup V}$  with
\begin{gather}\label{case3 par}
A=\frac{(\gamma-1)^2}{2},\qquad B=-\frac{(\gamma-\beta+n)^2}{2},\qquad C=k_1(\beta+n),\qquad D=-\frac{k_1^2}{2}.\end{gather}
\end{enumerate}
\end{theorem}

\begin{remark}
The parameters (\ref{case1 par}) are invariant under
$\beta\rightarrow 2-\beta$, $\gamma\rightarrow \gamma+1-\beta $;
compare with the parameters in the weight (\ref{weight sh}).
\end{remark}

\begin{remark}
Cases~2 and~3 in Theorem~\ref{thm: P5} follow from case 1 by
considering the B\"acklund transformation of Theorem \ref{thm:
BTr2}. Indeed, the compositions of the transformations
\[T_{1,1,-1}\circ T_{1,-1,1}=T_{1,1,1}\circ T_{1,-1,-1}\]
give the transformation
\[  Y_1 = y-\frac{2(\gamma-1)(y-1)^2 y}{t
y'-(1+\beta-2\gamma)y^2+(1+n-k_1t+\beta-2\gamma)y-n },
\] where
$y=y(t)$ solves {\rm P}$_{\textup V}$  with parameters (\ref{case1 par}) and $Y_1$
solves {\rm P}$_{\textup V}$  with (\ref{case2 par}). Similarly, the compositions
of the transformations \[
T_{1,1,-1}\circ
T_{-1,-1,1}=T_{1,1,1}\circ T_{-1,-1,-1}
\] give
\[ Y_2=w+\frac{2(\beta-\gamma)(y-1)^2y}{ty'-(1+\beta-2\gamma)y^2+(1+n-k_1
t+\beta-2\gamma )y-n},\] where $y=y(t)$ solves {\rm P}$_{\textup V}$  with
(\ref{case1 par}) and $Y_2$ solves {\rm P}$_{\textup V}$  with~(\ref{case3 par}).
 \end{remark}

Next we show that equation (\ref{nonlinear discrete vn}) can, in
fact, be obtained from the B\"acklund transformation of {\rm P}$_{\textup V}$  in
Theorem~\ref{thm: BTr2}. Let us take, for instance, the parameters
(\ref{case1 par}) and $k_1=1$ for simplicity (hence, $c=t$).
Suppose $y=y_n(t)$ solves {\rm P}$_{\textup V}$  with (\ref{case1 par}). Then by
considering the transformations $T_{1,-1,-1}\circ T_{-1,-1,1}\circ
T_{-1,-1,1}$ and $T_{1,-1,-1}\circ T_{1,1,-1}\circ T_{1,1,1}$ we
get new solutions of {\rm P}$_{\textup V}$   with parameters (\ref{case1 par}) for
$n+1$ and $n-1$ respectively. In particular,
\begin{gather*}
 y_{n+1} = 1-\frac{2t(n+\gamma)y}{(\beta-1)(t
y'+y(1+n+t-\beta+(\beta-1)y)-n)} \\
\phantom{y_{n+1} =}{}  +
\frac{2t(1+n-\beta+\gamma)y}{(\beta-1)(ty'+y(n-1+t+\beta+(1-\beta)
y )-n)}\end{gather*} and
\begin{gather*}
 y_{n-1} = 1+\frac{2t(\gamma+n-1)y}{(\beta-1)(ty'-y(1+n+t-\beta+(\beta-1)y)+n)}
\\
\phantom{y_{n-1} =}{} -  \frac{2t(n-\beta+\gamma)y}{(\beta-1)(ty'-y(n-1+t+\beta+(1-\beta)y)+n)}.
\end{gather*}
Expressing $v_{n\pm 1}$ in terms of $y_n$ by using~(\ref{case1}),
we can compute that (\ref{nonlinear discrete vn}) is identically
zero. Hence, the discrete system
(\ref{discrete1}), (\ref{discrete2}) can be obtained from the
B\"acklund transformation of {\rm P}$_{\textup V}$.

\subsection{Initial conditions}

In this section we study the initial conditions (\ref{initial}),
(\ref{initial sh}), (\ref{initial bi}) for the generalized Meixner
weight~(\ref{weight}) on the lattice~$\mathbb{N}$, on the shifted lattice
and on the bi-lattice respectively.

Let $n=0$. Since $\gamma\neq 1$, we get that $u_0=0$ from
$a_0^2=0.$ We also put $k_1=1$ and $c=t$. From~(\ref{*}) we have
that $v=v_0(t)$ satisf\/ies the f\/irst order nonlinear equation
\begin{gather}\label{Ric v}
t^2 v'=(\gamma-1)v^2+t(2-t+\beta-2\gamma)v+(\gamma-\beta)t^2.
\end{gather}
Since $b_0=\gamma-\beta+t-(\gamma-1)v_0/t,$ we can f\/ind $v_0$ for
(\ref{initial}), (\ref{initial sh}) and (\ref{initial bi}). We can
verify that all of them satisfy (\ref{Ric v}) using formulas for
the conf\/luent hypergeometric functions from~\cite{AbramowitzStegun}. Note that~(\ref{initial bi}) depends on
an arbitrary parameter~$\tau$.

The f\/ifth Painlev\'e equation (\ref{P5}) with parameters
(\ref{case2 par}) with $k_1=1$ and $n=0$ has particular solutions
which solve the following f\/irst order nonlinear equation:
\begin{gather}\label{Ric w}
t y'=(\beta-\gamma)y^2+(t-1-\beta+2\gamma)y+1-\gamma.
\end{gather}
Substituting expression (\ref{case2}) in (\ref{Ric v}) and
assuming that $y$ satisf\/ies (\ref{Ric w}), we indeed f\/ind that the
equation is satisf\/ied. We also f\/ind that $v(t)=t/y(t)$. Thus, the
initial conditions~(\ref{initial}),~(\ref{initial sh}),~(\ref{initial bi}) for the generalized Meixner weight
(\ref{weight}) on the lattice $\mathbb{N}$, on the shifted lattice and on
the bi-lattice respectively are related to solutions of the f\/irst
order dif\/ferential equation (\ref{Ric w}), which, in turn,
satisf\/ies {\rm P}$_{\textup V}$.

\section{A remark on the B\"acklund transformation}

 In this section we study the B\"acklund
transformation of the f\/ifth Painlev\'e equation (P$_{\rm V}$) and
f\/ind when a linear combination of two solutions is also a solution
of P$_{\rm V}$. In particular, we show that if $y_1$ and $y_2$
are solutions of P$_{\rm V}$
 obtained from a solution $y$ by certain
B\"acklund transformations, then there is a constant $M\neq 0, 1$
such that $M  y_1+(1-M) y_2$ is also a solution of P$_{\rm
V}$.

\begin{example}
In \cite{LiesGFWalter} it is shown that if $y:=y_n(t)$ (related to
the  recurrence coef\/f\/icients of the generalized Meixner
polynomials) is the solution of~(\ref{P5}) with
\begin{gather*}\label{first case}
A=\frac{(\beta-1)^2}{2},\qquad B=-\frac{(\beta+n)^2}{2},\qquad C=2n,\qquad D=-2,
\end{gather*}
then
one can show that  $y_{n+1}=y_{n+1}(t) $ given by
\begin{gather}\nonumber
 y_{n+1}= 1+\frac{4(n+1)ty}{(\beta-1)(ty'+(2t+2\beta-1+n-(\beta-1)y)y-n-\beta)}
\\
\phantom{y_{n+1}=}{}  - \frac{4t(n+\beta)y}{(\beta-1)(t
y'+(1+n+2t+(\beta-1)y)y-n-\beta)}\label{example1}
\end{gather}
 is the solution of
(\ref{P5}) with
\[A=\frac{(\beta-1)^2}{2},\qquad B=-\frac{(\beta+n+1)^2}{2},\qquad C=2(n+1),\qquad D=-2.\]
 It can be observed that  the transformation~(\ref{example1})
can be written in the following form: \[y_{n+1}=M
 y_1+(1-M) y_2,\qquad M=\frac{n+1}{1-\beta},\] where
\[y_1=T_{1,-1,1}y=1+\frac{4t y
}{n+\beta-(n-1+2\beta+2t)y+(\beta-1)y^2-t y'}
\] is a solution of
(\ref{P5}) with
\[A_1=\frac{(n+1)^2}{2},\qquad B_1=-\frac{1}{2},\qquad C_1=-2(2\beta+n-1),\qquad D_1=-2\]
and \[y_2=T_{-1,-1,1}y=1+\frac{4t y
}{n+\beta-(1+n+2t)y-(\beta-1)y^2-t y'}\] is a solution of
(\ref{P5}) with
\[A_2=\frac{(n+\beta)^2}{2},\qquad B_2=-\frac{\beta^2}{2},\qquad C_2=-2(n+1),\qquad D_2=-2.\]
\end{example}

Similarly, if
\begin{gather*}
y_1=T_{1,1,-1}y=1+\frac{4t y
}{n+\beta-(1+n+2t)y-(\beta-1)y^2+t y '},\\
 y_2=T_{-1,1,-1}y=1+\frac{4t y
}{n+\beta-(n-1+2\beta+2t)y+(\beta-1)y^2+t y'}
\end{gather*}
 are solutions of
(\ref{P5}) with
\[A_1=\frac{(\beta+n-1)^2}{2},\qquad B_1=-\frac{(\beta-1)^2}{2},\qquad C_1=-2(n+1),\qquad D_1=-2\]
and
\[A_2=\frac{n^2}{2},\qquad B_2=0,\qquad C_2=-2(2\beta+n-1),\qquad D_2=-2,\]
respectively, then \[y_{n-1}=M y_1+(1-M) y_2,
\qquad M=\frac{(\beta+n-1)}{\beta-1}\] is a solution of (\ref{P5})
with
\[A=\frac{(\beta-1)^2}{2},\qquad B=-\frac{(\beta+n-1)^2}{2},\qquad C=2(n-1),\qquad D=-2.\]
Such observations motivate us to study the question when the sum
$M y_1+(1-M)y_2$ of two solutions of P$_{\rm V}$  is also a
solution of the same equation. Clearly, we impose the conditions that
$M\neq 0 $ and $M\neq 1.$

\begin{theorem}
Let $y=y(t)$ be a solution of P$_{\rm V}$  with parameters
$A$, $B$, $C$, $D=-2$ and
\[y_1=T_{\varepsilon_1,\varepsilon_2,\varepsilon_3}y,\qquad y_2=T_{\delta_1,\delta_2,\delta_3}y,\]
where $\varepsilon_j^2=\delta_j^2=1$ and
$\varepsilon_j\neq\delta_j$ for some $j\in\{1,2,3\}.$ Then
\[v=M\,y_1+(1-M)\,y_2\] with $M\neq 0;1 $  is a solution of
{\rm P}$_{\rm V}$  with parameters $A_v$, $B_v$, $C_v$, $D_v=-2$ in the
following cases:
\begin{enumerate}\itemsep=0pt
\item[$1.$]
$\delta_1=\varepsilon_1$, $\delta_2=-\varepsilon_2$, $\delta_3=\varepsilon_3$
and
\begin{gather*}
M=\frac{2\varepsilon_1
\sqrt{2A}+2\varepsilon_2\sqrt{-2B}-\varepsilon_3C-2}{4\varepsilon_2\sqrt{-2B}},\\
 A_v=-B,\qquad B_v=\frac{2\varepsilon_1\sqrt{2A}-2A-1}{2},\qquad C_v=-C-2\varepsilon_3;
 \end{gather*}
\item[$2.$] $\delta_1=-\varepsilon_1$, $\delta_2=\varepsilon_2$, $\delta_3=\varepsilon_3$
and
\begin{gather*}
M=\frac{2\varepsilon_1\sqrt{2A}+2\varepsilon_2\sqrt{-2B}-\varepsilon_3C-2}{4\varepsilon_1\sqrt{2A}},\\
 A_v=A,\qquad B_v=\frac{2\varepsilon_2\sqrt{-2B}+2B-1}{2},\qquad C_v=C+2\varepsilon_3.
 \end{gather*}
\end{enumerate}

\end{theorem}

\begin{proof} The proof of this result is computational. We f\/irst
obtain that in case $\delta_3=-\varepsilon_3$ one gets only cases
$M=0$ and $M=1.$ In case $\delta_3=\varepsilon_3$ one needs to
consider 3 cases separately (since $\delta_1=\varepsilon_1,
\;\;\delta_2=\varepsilon_2$ gives a trivial result for the
function $v$): $\delta_1=\varepsilon_1$,
$\delta_2=-\varepsilon_2$; $\delta_1=-\varepsilon_1$,
$\delta_2=\varepsilon_2$; $\delta_1=-\varepsilon_1$,
$\delta_2=-\varepsilon_2$. However, in the last case we get
$M=0$ or $M=1.$
\end{proof}

The examples at the  beginning of the section correspond to the second
case of the theorem.
Similarly, we can get expressions for $y_{n+1}$ and $y_{n-1}$ in
the previous section by using this theorem.

Various properties of the repeated application of the B\"acklund
transformations are studied in~\cite{Gr1, Gr2, Gr3}. Repeated
applications of the B\"acklund transformations to the seed
solutions usually lead to  very cumbersome formulas. However, as
shown in this section, we can get linear dependence between three
solutions. Moreover, our formulas suggest that the function~$v$
has the same poles as~$y_1$ and~$y_2$ and, thus, they can be
useful to study various properties of the solutions. Other
Painlev\'e equations might have   similar properties so one can
try to study
when  for instance a linear combination of several solutions or a
product or a cross-ratio of several solutions is also a solution
(see also the representation of solutions in~\cite{Noumi}).
Although computationally   dif\/f\/icult, this deserves further study.

\subsection*{Acknowledgements}
GF is
 partially supported by Polish MNiSzW Grant N N201 397937.  WVA is supported by Belgian Interuniversity Attraction Pole P6/02, FWO
grant G.0427.09 and K.U.\ Leuven Research Grant OT/08/033.

\pdfbookmark[1]{References}{ref}
\LastPageEnding


\begin{thebibliography}{99}

\footnotesize\itemsep=0pt


\bibitem{AbramowitzStegun}
Abramowitz M., Stegun I.,
 Handbook of mathematical functions, Dover Publications, New York, 1965.

\bibitem{LiesGFWalter}
Boelen L.,  Filipuk G., Van Assche W.,
Recurrence coef\/f\/icients of generalized Meixner polynomials and Painlev\'e equations,
\href{http://dx.doi.org/10.1088/1751-8113/44/3/035202}{{\it J.~Phys.~A: Math. Theor.}} {\bf 44} (2011), 035202, 19~pages.


\bibitem{Chihara}
Chihara T.S.,
An introduction to orthogonal polynomials,
{\it Mathematics and its Applications}, Vol.~13, Gordon and Breach Science Publishers, New York~-- London~-- Paris, 1978.

\bibitem{GFWalterLaguerre}
Filipuk G.,  Van Assche W., Zhang L.,
On the recurrence coef\/f\/icients of semiclassical Laguerre polynomials,
\href{http://arxiv.org/abs/1105.5229}{arXiv:1105.5229}.

\bibitem{GFWalterCharlier}
Filipuk  G.,  Van Assche W.,
 Recurrence coef\/f\/icients of the generalized Charlier polynomials and the f\/ifth Painlev\'e equation,
\href{http://arxiv.org/abs/1106.2959}{arXiv:1106.2959}.

\bibitem{Gr1}
Gromak V., Filipuk  G.,
On functional relations between solutions of the f\/ifth Painlev\'e equation,
\href{http://dx.doi.org/10.1023/A:1019204329101}{{\it  Differ. Equ.}} {\bf 37} (2001), 614--620.

\bibitem{Gr2}
Gromak  V., Filipuk G.,
The B\"acklund transformations of the f\/ifth Painlev\'e equation and their applications,
\href{http://dx.doi.org/10.1080/13926292.2001.9637161}{{\it Math. Model. Anal.}} {\bf 6} (2001), 221--230.

\bibitem{Gr3}
Gromak V., Filipuk  G.,
B\"acklund transformations of the f\/ifth Painlev\'e equation and their applications,
in  Proceedings of the Summer School ``Complex Dif\/ferential and Functional Equations''  (Mekrij\"arvi, 2000),
{\it Univ. Joensuu Dept. Math. Rep. Ser.}, Vol.~5, Univ. Joensuu, Joensuu, 2003, 9--20.


\bibitem{GrLSh}
Gromak V.I., Laine  I.,  Shimomura S.,
 Painlev\'e  dif\/ferential equations in the complex plane,
{\it de Gruyter Studies in Mathematics}, Vol.~28, Walter de Gruyter \& Co., Berlin, 2002.

\bibitem{Magnus}
Magnus  A.P.,
Painlev\'e type dif\/ferential equations for the recurrence coef\/f\/icients of semi-classical orthogonal polynomials,
\href{http://dx.doi.org/10.1016/0377-0427(93)E0247-J}{{\it J. Comput. Appl. Math.}} {\bf 57} (1995), 215--237,
\href{http://arxiv.org/abs/math.CA/9307218}{math.CA/9307218}.

\bibitem{Noumi}
Noumi M.,
Painlev\'e equations through symmetry,
{\it Translations of Mathematical Monographs}, Vol.~223, American Mathematical Society, Providence, RI, 2004.

\bibitem{bilattice}
Smet C., Van Assche W.,
Orthogonal polynomials on a bi-lattice,  {\it Constr. Approx.},    to appear,
\href{http://arxiv.org/abs/1101.1817}{arXiv:1101.1817}.

\end{thebibliography}
\end{document}